\documentclass[12pt]{amsart}

\usepackage[left=3.5cm,top=3.5cm,right=3.5cm]{geometry}
\usepackage[T1]{fontenc}

\usepackage{enumerate}
\usepackage[shortlabels]{enumitem}
\setlist[enumerate]{label=(\alph*),ref=\alph*}

\usepackage{caption}
\usepackage{booktabs}

\usepackage{tikz}
\usetikzlibrary{cd,decorations.markings,arrows.meta,calc}

\usepackage{
  amsmath,
  amssymb,
  amsthm,
  amsfonts,
  mathtools,
  bbm,
  mathdots,
  xcolor,
  bm,
  mathrsfs,
  microtype,
  thmtools,
  graphicx,
  verbatim,
}

\allowdisplaybreaks

\usepackage[normalem]{ulem}

\definecolor{cadmiumgreen}{rgb}{0.0, 0.42, 0.24}
\usepackage[backref, hidelinks, colorlinks, citecolor=cadmiumgreen, hypertexnames=false]{hyperref}

\usepackage[alphabetic, msc-links, nobysame, lite, backrefs]{amsrefs}

\usepackage[capitalize]{cleveref}

\newtheorem{thm}{Theorem}[section]
\newtheorem{prop}[thm]{Proposition}
\newtheorem{lem}[thm]{Lemma}

\theoremstyle{definition}
\newtheorem{definition}[thm]{Definition}
\newtheorem{remark}[thm]{Remark}
\newtheorem{ex}[thm]{Example}

\crefformat{section}{\S{#2#1#3}}
\Crefformat{section}{\S{#2#1#3}}
\crefrangeformat{section}{\S\S{{#3#1#4}} to #5#2#6}
\Crefrangeformat{section}{\S\S{{#3#1#4}} to #5#2#6}
\crefmultiformat{section}{\S{\S{#2#1#3}}}{ and~#2#1#3}{, #2#1#3}{, and~#2#1#3}
\Crefmultiformat{section}{\S{\S{#2#1#3}}}{ and~#2#1#3}{, #2#1#3}{, and~#2#1#3}
\crefformat{equation}{(#2#1#3)}
\Crefformat{equation}{(#2#1#3)}
\crefrangeformat{equation}{(#3#1#4) to (#5#2#6)}
\Crefrangeformat{equation}{(#3#1#4) to (#5#2#6)}
\crefmultiformat{equation}{(#2#1#3)}{ and~(#2#1#3)}{, (#2#1#3)}{, and~(#2#1#3)}
\Crefmultiformat{equation}{(#2#1#3)}{ and~(#2#1#3)}{, (#2#1#3)}{, and~(#2#1#3)}
\crefname{ex}{Example}{Examples}

\usepackage{style}

\newcommand{\pL}{\mathcal{L}}
\newcommand{\Z}{\mathbb{Z}}

\newcommand{\R}{\mathbb{R}}

\newcommand{\A}{\mathbb{A}}
\newcommand{\Gm}{\mathbb{G}_m}

\newcommand{\F}{\mathbb{F}}
\newcommand{\K}{\mathbb{K}}
\newcommand{\I}{\mathrm{I}_{1,6}}
\newcommand{\ip}[2]{\langle #1,#2\rangle}
\DeclareMathOperator{\conv}{conv}
\DeclareMathOperator{\Trop}{Trop}
\DeclareMathOperator{\init}{in}

\title[Minuscule Coxeter matroids and the $E_6$-Grassmannian]{The $E_6$-Grassmannian, its minuscule Coxeter matroids, and the lines on a cubic surface}
\author{Daniel Corey}
\address{Department of Mathematics, Embry-Riddle Aeronautical University, Daytona Beach, FL, USA}
\email{daniel.corey@erau.edu}
\thanks{Partially supported by the AMS-Simons Research Enhancement Grant for PUI Faculty.}
\date{\today}

\begin{document}

\begin{abstract}
We study the minuscule Coxeter matroids of type $E_6$, whose ground set is the set of $27$ lines on a smooth cubic surface. We characterize them by a tritangent exchange property and classify them into $230$ orbits under the Weyl group. Every nonempty matroidal stratum of the $E_6$-Grassmannian is an open subscheme of a torus, and all orbits but one are realizable over every infinite field. Finally, we compute the tropical $E_6$-Grassmannian and its $60$ classes of toric degenerations. This is an extended abstract for a forthcoming paper.
\end{abstract}

\maketitle

\section{Introduction}
\label{sec:intro}

The $E_6$-Grassmannian is a remarkable variety that can be approached from several directions: through the exceptional Jordan algebra over the octonions (where it appears as the projective plane over the complexified octonions and is therefore also called the Cayley plane), through the representation theory of $E_6$, and through the $27$ lines on a smooth cubic surface. 
We show that the last of these gives an illuminating way to understand the $(E_6,P_6)$--matroids and their strata in the $E_6$-Grassmannian. 
Here $P_6$ is the maximal parabolic subgroup associated with the sixth simple root, numbered as in Bourbaki~\cite[Plate~V]{Bourbaki}. 
The matroid axiom becomes an exchange property for tritangent planes, the defining quadrics of $E_6/P_6$ are signed sums over the tritangent planes through a line, and the analysis of each stratum is governed by how the tritangent planes through different lines meet.

The homogeneous coordinates of $E_6/P_6\subseteq\mathbb P^{26}$ are indexed by the $27$ weights of the minuscule representation of $E_6$, whose convex hull is the Gosset polytope $2_{21}$; the $(E_6,P_6)$--matroids are certain subsets of its vertices. 
Coxeter identified the vertices of $2_{21}$ with the $27$ lines on a smooth cubic surface~\cite{Coxeter1940}, and Lee realized $2_{21}$ as the convex hull of the classes of those lines in the Picard lattice~\cite{Lee}. 
We start from Lee's description and project it to the orthogonal complement of the canonical class, which places $2_{21}$ in the weight lattice $M$ of $E_6$. 
The inner product on $M$ then turns combinatorial statements about $2_{21}$ into statements about lines on a cubic surface, and this dictionary drives all of our results.

\subsection{The combinatorics of \texorpdfstring{$(E_6,P_6)$}{(E6,P6)}--matroids}

Using the inner product on $M$, we characterize $(E_6,P_6)$--matroids by a \emph{tritangent exchange property} (\cref{thm:char}): a set $Q$ of lines is a matroid if and only if, whenever two lines of $Q$ meet, the third line of their tritangent plane lies on a second tritangent plane whose other two lines also belong to $Q$. 
By the Gelfand--Serganova theorem, a set of lines is a matroid exactly when every edge of its convex hull is parallel to a root of $E_6$; testing this condition on every orbit of subsets of the $27$ lines shows that there are $1\,335\,373$ $(E_6,P_6)$--matroids, forming $230$ orbits under $W(E_6)$ (\cref{thm:class}).

\subsection{Matroidal strata of the \texorpdfstring{$E_6$}{E6}-Grassmannian}

As for $\mathrm{Gr}(k,n)$, the support of every point of $E_6/P_6$ is an $(E_6,P_6)$--matroid, and the points with a given support $Q$ form a locally closed subscheme $\cM_Q$ defined over $\Z$, the \emph{matroidal stratum} of $Q$ (\cref{sec:strata}). 
In a chart $\{p_y\neq0\}$, the coordinates of the $16$ lines skew to $y$ are free, and the coordinates of the $10$ lines meeting $y$ are quadrics in them (\cref{thm:bigcell}); so each stratum is cut out of a torus by the vanishing of some of these quadrics. 
Our main geometric result (\cref{thm:main}) is that these strata are as simple as possible. 
Up to symmetry, exactly one $(E_6,P_6)$--matroid is realizable over no field. 
Every other stratum is isomorphic over $\Z$ to an open subscheme of an algebraic torus, so it is smooth, irreducible and rational, and the matroid is realizable over every infinite field. 
Both the non-realizable matroid and the one stratum that requires a separate argument are explained by a single linear relation among the chart quadrics.

\subsection{Tropicalization and toric degenerations}

Finally, we compute the tropicalization of the open stratum of $E_6/P_6$ (\cref{thm:trop}). 
It is a simplicial fan with $12\,885$ rays and $8\,712\,360$ maximal cones, in $10$ and $421$ symmetry classes, and every cone induces a subdivision of $2_{21}$ into matroid polytopes; the classification of $(E_6,P_6)$--matroids (\cref{thm:class}) supplies the starting cone that the traversal algorithm needs. 
Exactly $60$ classes of maximal cones are prime, and each gives a flat degeneration of $E_6/P_6$ to a projectively normal Gorenstein toric Fano variety whose polytope has $27$ lattice points and normalized volume $78$ (\cref{thm:toric}). 
One class is the Hibi degeneration~\cite{GonciuleaLakshmibai}, and two give the chain polytope of the minuscule $E_6$ poset.

\section{\texorpdfstring{Coxeter $(E_6,P_6)$}{(E6,P6)}--matroids}
\label{sec:matroids}

Let $\I=\Z h\oplus\Z a_1\oplus\cdots\oplus\Z a_6$ be the lattice with symmetric bilinear form $h\cdot h=1$, $a_i\cdot a_j=-\delta_{ij}$, $h\cdot a_i=0$, and put $K=-3h+\sum_i a_i$, so that $K\cdot K=3$. 
This is the Picard lattice of a smooth cubic surface $S\subseteq\mathbb P^3$ with its canonical class \cite[Ch.~9]{Dolgachev}. 
With this connection, we use classical vocabulary below (lines, tritangent planes, double sixes). 
The orthogonal complement $K^{\perp}\subseteq\I$ is the $E_6$ root lattice, whose $72$ roots are the vectors $\alpha$ with $\alpha\cdot\alpha=-2$.

\begin{definition}
Let $\pi\colon\I\otimes\R\to K^{\perp}\otimes\R$ be the orthogonal projection, $\pi(x)=x-\tfrac13(x\cdot K)\,K$. 
Set
\[
	M=\pi(\I),\qquad M_\R=M\otimes\R=K^{\perp}\otimes\R,\qquad \ip{u}{v}=-\,u\cdot v ,
\]
so that $\ip{\cdot}{\cdot}$ is a positive definite $W(E_6)$-invariant inner product on $M_\R$.
\end{definition}

Since $\I$ is unimodular and $K$ is primitive, $M$ is the dual of $K^{\perp}$ with respect to $\ip{\cdot}{\cdot}$. 
Thus $M$ is the weight lattice of $E_6$; it contains the root lattice $K^{\perp}$ with index $3$, and $M/K^{\perp}\cong\Z/3\Z$ is generated by the class of $\pi(a_1)=a_1+\tfrac13K$. 
The Weyl group $W(E_6)$, generated by the reflections in the roots, acts on $M$.
The $27$ vectors $x\in\I$ with $x\cdot x=-1$ and $x\cdot K=-1$ are the classes of the $27$ lines on the cubic surface,
\[
	a_i,\qquad b_i = 2h-\textstyle\sum_{j\neq i}a_j,\qquad c_{ij}=h-a_i-a_j ,
\]
and each projects to $\pi(x)=x+\tfrac13K$. We keep the names $a_i,b_i,c_{ij}$ for these projections and write $\pL\subseteq M$ for the resulting set of $27$ vectors. 
Equivalently,
\[
	\pL=\bigl\{u\in \pi(a_1)+K^{\perp} : \ip{u}{u}=\tfrac43\bigr\},
\]
the set of weights of the $27$-dimensional minuscule representation $V$ of $E_6$ (the other nontrivial coset of $K^\perp$ gives the weights of $V^{*}$).
Following~\cite{Lee}, the Gosset polytope is
\[
	2_{21}=\conv(\pL)\subseteq M_\R ,
\]
a lattice polytope in $M$ centered at the origin, with $W(E_6)$ as its full symmetry group. 
From $\ip{\pi(x)}{\pi(y)}=-x\cdot y+\tfrac13$ we read off that two distinct elements $x,y\in\pL$ satisfy $\ip{x}{y}\in\{\tfrac13,-\tfrac23\}$, according as the corresponding lines are skew or meet, and that $x,y$ span an edge of $2_{21}$ exactly when $\ip{x}{y}=\tfrac13$, that is, when $x-y$ is a root.

Fix $y\in\pL$. 
The inner product partitions the remaining elements as
\begin{equation}
\label{eq:partition}
	\pL\setminus\{y\} = N(y)\sqcup\lozenge(y),\quad
	N(y)=\{x : \ip{x}{y}=\tfrac13\}, \quad \lozenge(y)=\{x: \ip{x}{y}=-\tfrac23\},
\end{equation}
where $N(y)$ is the set of $16$ neighbors of $y$ in $2_{21}$, and the $10$ elements of $\lozenge(y)$ are the vertices of the cross-polytope facet of $2_{21}$ opposite to $y$. 
Two vertices $x,x'$ of $\lozenge(y)$ are antipodal in that cross-polytope if and only if $\ip{x}{x'}=-\tfrac23$, if and only if
\[
	x+x'+y=0 \quad\text{in } M.
\]
More generally, a triple $\{x,y,z\}$ of distinct elements of $\pL$ is called a \textit{tritangent plane} if $x+y+z = 0$ in $M$. 
In the cubic surface description, $-K$ is the class of a hyperplane section of $S\subseteq\mathbb P^3$, so $x+y+z=0$ in $M$ (equivalently, the classes sum to $-K$ in $\I$) means that the three lines are coplanar: some plane $H\subseteq\mathbb P^3$ cuts out $H\cap S=x\cup y\cup z$. 

So the five antipodal pairs of $\lozenge(y)$ are the five tritangent planes through $y$, and there are $45$ tritangent planes in all. 
We write $x^{*}=-x-y$ for the antipode of $x$ in $\lozenge(y)$.

\begin{ex}
\label{ex:a6}
For $y=a_6$ we get $\lozenge(a_6)=\{b_1,\dots,b_5,c_{16},\dots,c_{56}\}$ with antipodal pairs $\{b_i,c_{i6}\}$, and $N(a_6)$ consists of $a_1,\dots,a_5,b_6$ and the ten $c_{ij}$ with $i,j\leq 5$.
\end{ex}

By the Gelfand--Serganova theorem~\cite{BorovikGelfandWhite}, an $(E_6,P_6)$--\emph{matroid} is a nonempty subset $Q\subseteq\pL$ such that every edge of the polytope $\conv(Q)$ is parallel to a root of $E_6$. 
It is a \emph{strong} $(E_6,P_6)$--matroid if for all distinct $u,v\in Q$ there is a root $\alpha$ whose reflecting hyperplane separates $u$ and $v$ and such that the reflections of $u$ and $v$ in $\alpha$ both lie in $Q$.

\begin{definition}
\label{def:tritangent}
A nonempty $Q\subseteq\pL$ has the \emph{tritangent exchange property} if, for every $a,b\in Q$ with $\ip{a}{b}=-\tfrac23$, there are $c,d\in Q$ with $\ip{c}{d}=-\tfrac23$, $\{c,d\}\neq\{a,b\}$, and
\[
	a+b=c+d \quad\text{in } M.
\]
\end{definition}

Geometrically, if $a$ and $b$ are two meeting lines of $Q$ and $y=-a-b$ is the third line of their tritangent plane, the property asks for a second tritangent plane $\{c,d,y\}$ through $y$ with $c,d\in Q$; in other words, a tritangent plane through $y$ with two lines in $Q$ can always be exchanged for another one. 
Since $\pi$ is injective on each level set $\{x\cdot K=k\}$ of $\I$, the relation $a+b=c+d$ may be verified among the corresponding classes in $\I$.

\begin{thm}
\label{thm:char}
For a nonempty $Q\subseteq\pL$ the following are equivalent.
\begin{enumerate}[label=(\roman*)]
	\item $Q$ is an $(E_6,P_6)$--matroid.
	\item $Q$ is a strong $(E_6,P_6)$--matroid.
	\item For every $y\in\pL$, the number of antipodal pairs of $\lozenge(y)$
	contained in $Q$ is not equal to one. 
	\item $Q$ has the tritangent exchange property.
\end{enumerate}
\end{thm}

\begin{proof}[Indication of proof]
The equivalence of (i) and (ii) is
\cite[Thm.~4.22]{GrossKuhnLuber}, and the equivalence of (ii) and (iii) is \cite[Prop.~3.15]{CalvertDermenjianFinkSmith}. For a direct proof of (i)$\Rightarrow$(iii), if $\lozenge(y)\cap Q$ contains exactly one antipodal pair $\{x,x^*\}$, an explicit linear functional on the cross-polytope shows that $[x,x^*]$ is an edge of $\conv(\lozenge(y)\cap Q)$, which is a face of $\conv(Q)$ because $\lozenge(y)$ spans a facet; but $\ip{x-x^*}{x-x^*}=4$, so $x-x^*$ is not parallel to a root. 
For (iii)$\Leftrightarrow$(iv), if $\ip{a}{b}=-\frac23$ then $y=-a-b$ lies in $\pL$ by the description of $\pL$ as the vectors of norm $\tfrac43$ in a coset, and $\{a,b\}$ is an antipodal pair of $\lozenge(y)$; conversely $\ip{c+d}{y}=\ip{a+b}{y}=-\tfrac43$ forces $c,d\in\lozenge(y)$.
\end{proof}

\medskip
\noindent In particular the adjective ``strong'' is redundant for $(W(E_6),P_6)$ and for all minuscule types by \cite{GrossKuhnLuber}, but it fails for general pairs $(W,P)$.

\begin{thm}
\label{thm:class}
There are exactly $1\,335\,373$ $(E_6,P_6)$--matroids, forming $230$ orbits under $W(E_6)$. 
By the dimension of $\conv(Q)$ the orbits are distributed as follows.
\begin{center}
\begin{tabular}{c|ccccccc}
	\toprule
	$\dim\conv(Q)$ & 0 & 1 & 2 & 3 & 4 & 5 & 6\\
	\midrule
	$W(E_6)$-orbits & 1 & 1 & 2 & 4 & 11 & 30 & 181\\
	\bottomrule
\end{tabular}
\end{center}
\end{thm}

\begin{proof}[Indication of proof]
Order $\pL$ and represent each orbit of subsets by its lexicographically minimal element. Every orbit of
$(k+1)$-subsets contains $R\cup\{v\}$ for some minimal representative $R$ of size $k$, so the representatives can be grown one line at a time; those of size at least $14$ are complements of those of size at most $13$. 
This yields the $5\,486$ orbits of subsets of $\pL$, and for each representative we test whether every edge of its convex hull is parallel to a root.
\end{proof}

\medskip
Write $\Gamma_Q$ for the graph on $Q$ joining $x$ and $z$ when $\ip{x}{z}=-\tfrac23$. 
Each edge $\{x,z\}$ of $\Gamma_Q$ determines the line $-x-z$, the third line of their tritangent plane. 
One matroid plays a special role below:
\[
	Q_1=\{a_1,\dots,a_6,\ b_1,\dots,b_6,\ c_{12},c_{13},c_{23},c_{45},c_{46},c_{56}\}.
\]
Its complement is the union of three disjoint tritangent planes, namely
$\{c_{14},c_{25},c_{36}\}$, $\{c_{15},c_{26},c_{34}\}$ and
$\{c_{16},c_{24},c_{35}\}$, and every vertex of $\Gamma_{Q_1}$ has degree $7$.

\section{Matroidal strata of the \texorpdfstring{$E_6$}{E6}-Grassmannian}
\label{sec:strata}

\subsection{Equations and charts}
\label{sec:equations}

Let $R=\Z[p_\ell:\ell\in\pL]$ and let $C$ be the skew-symmetric $6\times6$ matrix with entries $p_{c_{ij}}$ above the diagonal. 
The $E_6$-invariant cubic on the $27$-dimensional minuscule representation, which goes back to Cartan's thesis~\cite{Cartan} (see~\cite{ManivelLines}, and also \cite[Prop. 4.15]{CalvertDermenjianFinkSmith}), is a signed sum of $45$ monomials $p_xp_yp_z$, one for each tritangent plane, that is, each triple with $x+y+z=0$ in $M$~\cite{Faulkner,ManivelLines}. With $\mathbf a=(p_{a_1},\dots,p_{a_6})^T$ and $\mathbf b=(p_{b_1},\dots,p_{b_6})^T$ it is
\[
	D=\mathbf{a}^{T}C\,\mathbf{b}+\mathrm{Pf}(C).
\]
  The $30$ monomials of $\mathbf a^TC\,\mathbf b$ are the tritangent planes $\{a_i,b_j,c_{ij}\}$, and the $15$ monomials of $\mathrm{Pf}(C)$ are the tritangent planes $\{c_{ij},c_{kl},c_{mn}\}$ with $\{i,j\}\sqcup\{k,l\}\sqcup\{m,n\}=\{1,\dots,6\}$.
The $E_6$-Grassmannian is the singular locus of the hypersurface defined by $D=0$~\cite{IlievManivel}, so it is cut out by the partial derivatives of $D$:
\[
	E_6/P_6=\mathrm{Proj}\,R/\langle D_y : y\in\pL\rangle,
	\qquad D_y=\frac{\partial D}{\partial p_y}.
\]
Differentiating the monomials shows that each quadric is a signed sum over the tritangent planes through $y$:
\begin{equation}
\label{eq:Dy}
	D_y=\sum_{\{x,x^{*}\}\subseteq\lozenge(y)}\pm\,p_x\,p_{x^{*}} .
\end{equation}
(In the full version we match these coordinates with the model of $E_6/P_6$ inside the exceptional Jordan algebra over the octonions, under which $D$ becomes the determinant.)

Fix $y\in\pL$ and let $U_y=\{p_y\neq0\}\subseteq E_6/P_6$. 
For $v\in\lozenge(y)$ we have $y\in\lozenge(v^*)$, so by \eqref{eq:Dy} the quadric $D_{v^*}$ contains the monomial $p_vp_y$; write $D_{v^*}=\pm(p_vp_y-g_v)$. 
The other four tritangent pairs $\{u,u^*\}$ of $\lozenge(v^*)$ satisfy $u+u^*=v+y$, hence lie in $N(y)$, so $g_v\in\Z[p_u:u\in N(y)]$ is a signed sum of four products $p_up_{u^*}$.

\begin{thm}[Big cell]
\label{thm:bigcell}
For every $y\in\pL$ the morphism
\[
	\A^{16}=\operatorname{Spec}\Z[p_u : u\in N(y)]\ \longrightarrow\ U_y,
	\qquad p_y=1,\quad p_v=g_v\ \ (v\in\lozenge(y)),
\]
which keeps the coordinates $p_u$, $u\in N(y)$, is an isomorphism.
\end{thm}

\begin{proof}[Indication of proof]
On $U_y$, with $p_y=1$, the equations $D_{v^*}=0$ say exactly that $p_v=g_v$. 
It remains to check that the other $17$ quadrics vanish identically after this substitution: $D_y$ becomes a quartic $\sum\pm g_vg_{v^*}$, and each $D_w$, $w\in N(y)$, becomes a relation $\sum\pm p_ug_v$ that is linear in the $g_v$. 
This is a straightforward computation in one chart, and $W(E_6)$ permutes the charts transitively by signed permutations of the coordinates.
\end{proof}

\subsection{Supports and strata}

For a field $\K$ and a $\K$-point $p$ of $E_6/P_6$, let $\supp(p)=\{\ell\in\pL : p_\ell\neq0\}$.

\begin{prop}
The support of every point of $E_6/P_6$ is an $(E_6,P_6)$--matroid.
\end{prop}

\noindent Indeed, if $\lozenge(y)\cap\supp(p)$ contained exactly one tritangent pair, then by \eqref{eq:Dy} the quadric $D_y$ would evaluate to a single nonzero monomial at $p$. 
This is analogous to the fact that the support of a point of $\mathrm{Gr}(k,n)$ is a matroid, and it motivates the following definition.

\begin{definition}
The \emph{matroidal stratum} of an $(E_6,P_6)$--matroid $Q$ is the locally closed subscheme
\[
	\cM_Q=\{p\in E_6/P_6 : \supp(p)=Q\},
\]
defined over $\Z$. We say $Q$ is \emph{realizable over} $\K$ if $\cM_Q$ has a $\K$-point.
\end{definition}

Thus $E_6/P_6=\bigsqcup_Q\cM_Q$, with $W(E_6)$ permuting the strata, and the open stratum $\cM_{\pL}$ is the locus where all $27$ coordinates are nonzero. 
These are analogous to the thin Schubert cells in~\cite{GGMS} of the ordinary Grassmannian.

\subsection{The equations of a stratum}

Fix $y\in Q$. Every point of $\cM_Q$ lies in the chart $U_y$, and setting the coordinates of $N(y)\setminus Q$ to zero in \cref{thm:bigcell} gives
\begin{equation}
\label{eq:stratum}
	\cM_Q\cong\bigl\{p\in\Gm^{Q\cap N(y)} :\ g_v(p)=0\ (v\in E_y),\quad
	g_v(p)\neq0\ (v\in\lozenge(y)\cap Q)\bigr\},
\end{equation}
where we keep writing $g_v$ after setting the coordinates of $N(y)\setminus Q$ to zero, and $E_y$ is the set of $v\in\lozenge(y)\setminus Q$ for which $g_v$ is not the zero polynomial, so that $g_v=0$ is a genuine equation. 
Each $g_v$ is a signed sum of monomials $p_up_{u^*}$ over antipodal pairs $\{u,u^*\}$ of $\lozenge(v^*)$, which are disjoint, so no variable appears twice; and by \cref{thm:char}(iii) $g_v$ is never a monomial. 
Hence each equation can be solved for any one of its variables.
To organize the elimination, note that the monomial $p_up_{u^*}$ of $g_v$ comes from two meeting lines of $N(y)$ with $u+u^*=v+y$. So we define the \emph{exchange graph} $\Gamma_Q(y)$ as the subgraph of $\Gamma_Q$ on $Q\cap N(y)$ whose edges are the pairs $\{u,u'\}$ with $u+u'-y\notin Q$. Its edges are exactly the monomials of the equations, and the label $u+u'-y\in E_y$ partitions them into one matching per equation. 
A \emph{leaf order} is an ordering of the matchings such that each one contains a vertex of degree one in the union of itself and its successors. 
The key combinatorial input is the following.

\begin{lem}
\label{lem:onepoint}
Let $v,v'\in\lozenge(y)$ be distinct and let $\{u,u^*\}\subseteq N(y)$ be a tritangent pair of $\lozenge(v^*)$. 
Then $\{u,u^*\}$ meets $\lozenge(v'^*)$ in no vertex if $v'=v^*$, and in exactly one vertex otherwise.
\end{lem}

This is a one-line computation, $\ip{u+u^*}{v'^*}=\ip{v+y}{v'^*}$ (both pairs sum to $-v^*$ in $M$), and it implies that two matchings of $\Gamma_Q(y)$ never have the same vertex set, so every chart with at most two equations has a leaf order. 
Eliminating one variable per equation along a leaf order identifies $\cM_Q$ with an open subscheme of a torus over $\Z$.

\begin{ex}
\label{ex:leaf}
Let $Q=\{a_1,\dots,a_5,b_1,\dots,b_4,b_6,c_{12},c_{13},c_{24},c_{35},c_{36},c_{45},c_{46}\}$, the $(E_6,P_6)$--matroid that supplies the starting cone in the proof of
\cref{thm:trop}, and take $y=a_3$. Then $Q\cap\lozenge(y)$ consists of $b_1$, $b_2$, $b_4$, $b_6$, $c_{13}$, $c_{35}$, $c_{36}$, and $Q\cap N(y)$ of $a_1$, $a_2$, $a_4$, $a_5$, $b_3$, $c_{12}$, $c_{24}$, $c_{45}$, $c_{46}$.
\Cref{fig:leaf}(a) shows $\Gamma_Q$, with $y$ on top, $Q\cap\lozenge(y)$ in the middle row and $Q\cap N(y)$ in the bottom row. 
Each of the $13$ edges
$\{u,u'\}$ of $\Gamma_Q$ on $Q\cap N(y)$ is labeled by the line $u+u'-y$; for instance $c_{12}+c_{46}=2h-a_1-a_2-a_4-a_6=b_5+a_3$ in $\I$.
\Cref{fig:leaf}(b) shows these edges with their labels. 
Seven labels lie in $Q$ (dashed edges), and these edges give the coordinates
\begin{gather*}
	p_{b_1}=p_{a_1}p_{b_3},\quad p_{b_2}=p_{a_2}p_{b_3},\quad p_{b_4}=p_{a_4}p_{b_3},\quad 
	p_{b_6}=-p_{c_{12}}p_{c_{45}},\\ 
	p_{c_{13}}=-p_{a_2}p_{c_{12}},\quad
	p_{c_{35}}=-p_{a_4}p_{c_{45}},\quad p_{c_{36}}=-p_{a_4}p_{c_{46}},
\end{gather*}
which are automatically units. 
The remaining six (solid) edges form $\Gamma_Q(a_3)$, with matchings $M_{b_5}$ (blue), $M_{c_{23}}$ (orange) and $M_{c_{34}}$ (green), so $E_y=\{b_5,c_{23},c_{34}\}$ and the equations are
\begin{gather*}
	g_{b_5}=p_{a_5}p_{b_3}+p_{c_{12}}p_{c_{46}}=0,\qquad
	g_{c_{23}}=p_{a_1}p_{c_{12}}-p_{a_4}p_{c_{24}}=0,\\
	g_{c_{34}}=p_{a_5}p_{c_{45}}-p_{a_2}p_{c_{24}}=0.
\end{gather*}
With three equations \cref{lem:onepoint} does not apply directly, but $\Gamma_Q(a_3)$ is the disjoint union of the paths $a_2c_{24}a_4$,
$c_{46}c_{12}a_1$ and $b_3a_5c_{45}$, so each matching contains a leaf of the whole graph and every order is a leaf order. 
Solving each equation for a leaf that occurs in no other equation gives
\[
	p_{a_1}=\frac{p_{a_4}p_{c_{24}}}{p_{c_{12}}},\qquad
	p_{a_2}=\frac{p_{a_5}p_{c_{45}}}{p_{c_{24}}},\qquad
	p_{b_3}=-\frac{p_{c_{12}}p_{c_{46}}}{p_{a_5}}.
\]
Every coordinate is now a Laurent monomial in
$p_{a_4},p_{a_5},p_{c_{12}},p_{c_{24}},p_{c_{45}},p_{c_{46}}$, so
$\cM_Q\cong\Gm^6$.
\end{ex}

\begin{figure}[t]
\centering
\tikzset{
  ln/.style={circle,draw,thick,minimum size=6.4mm,inner sep=0pt,font=\small},
  qed/.style={thick,densely dashed,gray},
  lab/.style={font=\scriptsize,fill=white,inner sep=1pt},
  mb5/.style={line width=1.6pt,blue!70!black},
  mc23/.style={line width=1.6pt,orange!85!black},
  mc34/.style={line width=1.6pt,green!50!black},
}
\begin{tikzpicture}
  \node[ln] (a2) at (0.00,0) {$a_{2}$};
  \node[ln] (c24) at (1.30,0) {$c_{24}$};
  \node[ln] (a4) at (2.60,0) {$a_{4}$};
  \node[ln] (c46) at (3.90,0) {$c_{46}$};
  \node[ln] (c12) at (5.20,0) {$c_{12}$};
  \node[ln] (a1) at (6.50,0) {$a_{1}$};
  \node[ln] (b3) at (7.80,0) {$b_{3}$};
  \node[ln] (a5) at (9.10,0) {$a_{5}$};
  \node[ln] (c45) at (10.40,0) {$c_{45}$};
  \node[ln,fill=green!15] (b1) at (0.20,1.9) {$b_{1}$};
  \node[ln,fill=green!15] (b6) at (1.65,1.9) {$b_{6}$};
  \node[ln,fill=green!15] (b2) at (3.10,1.9) {$b_{2}$};
  \node[ln,fill=green!15] (b4) at (4.55,1.9) {$b_{4}$};
  \node[ln,fill=green!15] (c35) at (6.00,1.9) {$c_{35}$};
  \node[ln,fill=green!15] (c13) at (7.45,1.9) {$c_{13}$};
  \node[ln,fill=green!15] (c36) at (8.90,1.9) {$c_{36}$};
  \node[ln,fill=yellow!40] (a3) at (4.55,3.4) {$a_3$};
  \foreach \b in {b1,b6,b2,b4,c35,c13,c36} \draw (a3)--(\b);
  \begin{scope}[gray!70]
  \foreach \u/\w in {b1/a2,b1/a4,b1/a5,b1/c12,b2/a1,b2/a4,b2/a5,b2/c12,b2/c24,b4/a1,b4/a2,b4/a5,b4/c24,b4/c45,b4/c46,b6/a1,b6/a2,b6/a4,b6/a5,b6/c46,c13/a1,c13/b3,c13/c24,c13/c45,c13/c46,c35/a5,c35/b3,c35/c12,c35/c24,c35/c46,c36/b3,c36/c12,c36/c24,c36/c45} \draw (\u)--(\w);
  \end{scope}
  \draw (b1.north) .. controls +(0.5,0.75) and +(-0.5,0.75) .. (c13.north);
  \draw (b6.north) .. controls +(0.5,0.75) and +(-0.5,0.75) .. (c36.north);
  \draw[] (a1.south) .. controls +(0.31,-0.45) and +(-0.31,-0.45) .. (b3.south);
  \draw[] (c12.south) .. controls +(0.31,-0.45) and +(-0.31,-0.45) .. (a1.south);
  \draw[] (a2.south) .. controls +(1.09,-1.89) and +(-1.09,-1.89) .. (b3.south);
  \draw[] (a2.south) .. controls +(0.77,-1.36) and +(-0.77,-1.36) .. (c12.south);
  \draw[] (a2.south) .. controls +(0.31,-0.45) and +(-0.31,-0.45) .. (c24.south);
  \draw[] (a4.south) .. controls +(0.77,-1.36) and +(-0.77,-1.36) .. (b3.south);
  \draw[] (c24.south) .. controls +(0.31,-0.45) and +(-0.31,-0.45) .. (a4.south);
  \draw[] (a4.south) .. controls +(1.09,-1.89) and +(-1.09,-1.89) .. (c45.south);
  \draw[] (a4.south) .. controls +(0.31,-0.45) and +(-0.31,-0.45) .. (c46.south);
  \draw[] (b3.south) .. controls +(0.31,-0.45) and +(-0.31,-0.45) .. (a5.south);
  \draw[] (a5.south) .. controls +(0.31,-0.45) and +(-0.31,-0.45) .. (c45.south);
  \draw[] (c12.south) .. controls +(0.77,-1.36) and +(-0.77,-1.36) .. (c45.south);
  \draw[] (c46.south) .. controls +(0.31,-0.45) and +(-0.31,-0.45) .. (c12.south);
  \node[font=\small] at (4.55,-2.25) {(a) $\Gamma_Q$};
\end{tikzpicture}

\bigskip
\begin{tikzpicture}
  \node[ln] (a2) at (0.00,0) {$a_{2}$};
  \node[ln] (c24) at (1.30,0) {$c_{24}$};
  \node[ln] (a4) at (2.60,0) {$a_{4}$};
  \node[ln] (c46) at (3.90,0) {$c_{46}$};
  \node[ln] (c12) at (5.20,0) {$c_{12}$};
  \node[ln] (a1) at (6.50,0) {$a_{1}$};
  \node[ln] (b3) at (7.80,0) {$b_{3}$};
  \node[ln] (a5) at (9.10,0) {$a_{5}$};
  \node[ln] (c45) at (10.40,0) {$c_{45}$};
  \draw[qed] (a1.north) .. controls +(0.31,0.45) and +(-0.31,0.45) .. node[lab]{$b_{1}$} (b3.north);
  \draw[mc23] (c12.south) .. controls +(0.31,-0.45) and +(-0.31,-0.45) .. node[lab,text=orange!85!black]{$c_{23}$} (a1.south);
  \draw[qed] (a2.north) .. controls +(1.09,1.89) and +(-1.09,1.89) .. node[lab]{$b_{2}$} (b3.north);
  \draw[qed] (a2.north) .. controls +(0.77,1.36) and +(-0.77,1.36) .. node[lab]{$c_{13}$} (c12.north);
  \draw[mc34] (a2.south) .. controls +(0.31,-0.45) and +(-0.31,-0.45) .. node[lab,text=green!50!black]{$c_{34}$} (c24.south);
  \draw[qed] (a4.north) .. controls +(0.77,1.36) and +(-0.77,1.36) .. node[lab]{$b_{4}$} (b3.north);
  \draw[mc23] (c24.south) .. controls +(0.31,-0.45) and +(-0.31,-0.45) .. node[lab,text=orange!85!black]{$c_{23}$} (a4.south);
  \draw[qed] (a4.north) .. controls +(1.09,1.89) and +(-1.09,1.89) .. node[lab]{$c_{35}$} (c45.north);
  \draw[qed] (a4.north) .. controls +(0.31,0.45) and +(-0.31,0.45) .. node[lab]{$c_{36}$} (c46.north);
  \draw[mb5] (b3.south) .. controls +(0.31,-0.45) and +(-0.31,-0.45) .. node[lab,text=blue!70!black]{$b_{5}$} (a5.south);
  \draw[mc34] (a5.south) .. controls +(0.31,-0.45) and +(-0.31,-0.45) .. node[lab,text=green!50!black]{$c_{34}$} (c45.south);
  \draw[qed] (c12.north) .. controls +(0.77,1.36) and +(-0.77,1.36) .. node[lab]{$b_{6}$} (c45.north);
  \draw[mb5] (c46.south) .. controls +(0.31,-0.45) and +(-0.31,-0.45) .. node[lab,text=blue!70!black]{$b_{5}$} (c12.south);
  \node[font=\small] at (4.55,-1.4) {(b) $\Gamma_Q$ on $Q\cap N(y)$};
\end{tikzpicture}
\caption{The graphs of \cref{ex:leaf}.}
\label{fig:leaf}
\end{figure}

\begin{thm}
\label{thm:main}
Let $Q$ be an $(E_6,P_6)$--matroid. If $Q$ is $W(E_6)$-equivalent to
\[
	Q_{\mathrm{nr}}=Q_1\cup\{c_{14}\}
	=\{a_1,\dots,a_6,\ b_1,\dots,b_6,\ c_{12},c_{13},c_{14},c_{23},c_{45},c_{46},c_{56}\},
\]
then $\cM_Q$ is empty over every field. 
Otherwise $\cM_Q$ is isomorphic over $\Z$ to an open subscheme of an algebraic torus with nonempty fiber over every prime. 
In particular $Q$ is realizable over every infinite field, and $\cM_Q$ is smooth, irreducible and rational, of dimension $|Q\cap N(y)|-|E_y|$ for any
chart $y$ admitting a leaf order.
\end{thm}

\begin{proof}[Indication of proof]
A chart with at most two equations has a leaf order by \cref{lem:onepoint}. 
A check over the $230$ orbit representatives shows that every $Q$ other than $Q_1$ has a chart admitting a leaf order. 
In each chart of $Q_1$ the three equations are the $2\times2$ minors of a $2\times3$ matrix, which cut out a translate of an $8$-dimensional subtorus. 
For $Q_{\mathrm{nr}}$, a linear relation among the $g_v$ in the chart at $a_1$ forces $p_{a_4}g_{c_{14}}=0$, although both factors must be units. 
Finally, strong Gr\"obner bases over $\Z$ with leading coefficients $\pm1$ show that the fibers over $\operatorname{Spec}\Z$ are nonempty.
\end{proof}

\begin{center}
\begin{tabular}{c|ccccccccccccccccc}
	\toprule
	$\dim\cM_Q$ & 0&1&2&3&4&5&6&7&8&9&10&11&12&13&14&15&16\\
	\midrule
	orbits & 1&1&2&3&8&12&23&34&43&38&30&17&9&4&2&1&1\\
	\bottomrule
\end{tabular}
\captionof{table}{Dimensions of the nonempty matroidal strata, up to
$W(E_6)$.}
\label{tab:dims}
\end{center}

\begin{remark}
\label{rem:char2}
The $(E_6,P_6)$--matroids other than $Q_{\mathrm{nr}}$ might not be realizable over certain finite fields. For instance $\pL$ itself is not realizable over $\F_2$. 
Indeed, a point with all coordinates nonzero would have all coordinates equal to $1$, and then each quadric \eqref{eq:Dy} would equal a sum of five terms $\pm1$, which is odd. 
\end{remark}

\begin{remark}
\label{rem:frontier}
Although each stratum is smooth and irreducible, the strata do not form a stratification in the sense that the closure of a stratum is not necessarily the union of smaller strata. 
For example, for $Q=\{a_1,\dots,a_5,b_6,c_{12},\dots,c_{15},c_{23},c_{45}\}$ and $P=\{a_2,\dots,a_5,c_{12},\dots,c_{15}\}$ we have that $\overline{\cM_Q}$ meets $\cM_P$ without containing it.
\end{remark}

\section{Tropicalization and toric degenerations}
\label{sec:trop}

\subsection{The tropical \texorpdfstring{$E_6$}{E6}-Grassmannian}

Let $X^\circ=\cM_{\pL}$ be the open stratum, a closed subvariety of the dense torus $\Gm^{27}/\Gm$ of $\mathbb P^{26}$ which is a smooth and irreducible very affine variety.
Let $I=\langle D_y:y\in\pL\rangle$. 
The initial form $\init_w f$ picks out the terms of $f$ of minimal $w$-weight.  
Since $I$ is homogeneous, $\init_w I$ depends only on the class of $w$ in $\R^{27}/\R\mathbf 1\cong\R^{26}$, and $\Trop(X^\circ)\subseteq\R^{27}/\R\mathbf 1$ is the set of classes $w$ for which $\init_w I$ contains no monomial~\cite{MaclaganSturmfels}.
It is a pure $16$-dimensional fan. Its lineality space $\Lambda$ is the $6$-dimensional subspace given by the image of $M^*_\R$ under $\eta\mapsto(\eta(\ell))_{\ell\in\pL}$, that is, the space of affine-linear functions on $2_{21}$ modulo constants.

\begin{thm}
\label{thm:trop}
The restriction of the Gr\"obner fan of $I$ to  $\Trop(X^\circ)$ is simplicial and each of its maximal cones has multiplicity one. 
Modulo $\Lambda$ its $f$-vector is
\begin{align*}
	(&12\,885,\ 201\,261,\ 1\,471\,647,\ 6\,556\,740,\ 19\,592\,172,\\
	&40\,753\,818,\ 59\,005\,476,\ 57\,133\,440,\ 33\,247\,440,\ 8\,712\,360).
\end{align*}
Up to the action of $W(E_6)$ it is $(10,56,224,697,1639,2841,3558,3070,1639,421)$. 
For each $w$ in $\Trop(X^\circ)$, the regular subdivision of $2_{21}$ induced by $w$ is a subdivision into $(E_6,P_6)$--matroid polytopes.
\end{thm}

\begin{proof}[Indication of proof]
We computed the fan with
\texttt{gfan\_tropicaltraverse}~\cite{gfan}, up to the symmetry of $W(E_6)$.
This traversal needs a maximal cone to start from, and the heuristic for finding one, \texttt{gfan\_tropicalstartingcone}, did not work for $I$.
Instead, we started from the height $-\mathbf 1_Q$, the negative of the indicator vector of the matroid $Q$ of \cref{ex:leaf}. 
For the last statement, up to signs the quadrics $D_y$ are the strong exchange equations defining the Coxeter Dressian $\mathrm{Dr}(E_6,P_6)$, so $\Trop(X^\circ)\subseteq\mathrm{Dr}(E_6,P_6)$, and every point of the Dressian induces a matroidal subdivision by~\cite[Prop.~4.15]{CalvertDermenjianFinkSmith} and \cite[Thm.~4.22]{GrossKuhnLuber}. 
We also verified it directly on a representative of each of the $14\,155$ orbits of cones.
\end{proof}

\medskip
The ten orbits of rays have sizes $27$, $36$, $36$, $216$, $240$, $270$, $540$, $720$, $4\,320$ and $6\,480$; the two orbits of size $36$ are spanned, modulo $\Lambda$, by the indicator vectors of the $36$ double sixes and of their complements. 
The subdivisions induced by maximal cones have between $51$ and $69$ maximal cells.

\begin{prop}
\label{prop:dressian}
The inclusion $\Trop(X^\circ)\subseteq\mathrm{Dr}(E_6,P_6)$ is strict.
\end{prop}

\begin{proof}[Indication of proof]
Let $w$ be the vector with $w_\ell=-1$ for $\ell\in Q_{\mathrm{nr}}$ and $w_\ell=0$ otherwise. 
It induces a matroidal subdivision of $2_{21}$ into $27$ cells, one of which is $\conv(Q_{\mathrm{nr}})$, so $w\in\mathrm{Dr}(E_6,P_6)$. 
A direct calculation of the initial ideal verified that $w\notin\Trop(X^\circ)$.
\end{proof}

\subsection{Toric degenerations}
\label{sec:toric}

A maximal cone $\sigma$ of $\Trop(X^\circ)$ is \emph{prime} if $\init_w I$ is a prime ideal for $w$ in its relative interior. 
In that case the Gr\"obner degeneration $E_6/P_6$ to $V(\init_w I)$ is a flat degeneration to a variety which, having a dense torus orbit and a binomial ideal, is a projective toric variety $X_{A_\sigma}$; here $A_\sigma\subseteq\Z^{17}$ is the configuration of $27$ lattice points obtained by dividing $\Z^{27}$ by the lattice of exponent differences of the binomials. 
Let $\Delta_\sigma=\conv(A_\sigma)$.

\begin{thm}
\label{thm:toric}
Exactly $60$ of the $421$ $W(E_6)$-orbits of maximal cones of $\Trop(X^\circ)$ are prime. 
For each prime cone $\sigma$, the special fiber $X_{A_\sigma}$ is a projectively normal toric variety.
Consequently it is Gorenstein with $\omega_{X_{A_\sigma}}\cong\mathcal O(-12)$, and in particular Fano.
\end{thm}

\begin{proof}[Indication of proof]
For a representative $w$ of each maximal cone orbit we compute $\init_w I$ by a Gr\"obner basis computation; it is binomial for $388$ orbits, and exactly $60$ of these binomial ideals are prime. 
For these, \textsc{Oscar}~\cite{OSCAR} shows that $A_\sigma$ is normal, so $X_{A_\sigma}$ is arithmetically Cohen--Macaulay.
The Hilbert series of $X_{A_\sigma}$ is $(1+10t+28t^2+28t^3+10t^4+t^5)/(1-t)^{17}$. 
Its numerator is palindromic, so $X_{A_\sigma}$ is Gorenstein by \cite[Thm.~4.4]{StanleyHilbert}, and a toric computation gives $\omega_{X_{A_\sigma}}\cong\mathcal O(-12)$.
\end{proof}

\medskip
The $60$ polytopes share much of their structure. Each $\Delta_\sigma$ is a
$16$-dimensional lattice polytope whose only lattice points are its $27$ vertices, with $297$ edges; the $54$ missing pairs are the diagonals of the $27$ parallelograms recorded by the quadratic binomials of $\init_w I$. Its $h^*$-vector is $(1,10,28,28,10,1)$.  
The orbits differ in finer invariants. 
For $53$ of them $\init_w I$ is generated by the $27$ quadratic binomials, while the other $7$ need one, three or eight additional cubic binomials. 
The number of facets of $\Delta_\sigma$ ranges from $22$ to $38$ (\cref{tab:facets}), and the $60$ polytopes realize $53$ combinatorial types.

\begin{center}
\begin{tabular}{c|ccccccccccccc}
	\toprule
	facets of $\Delta_\sigma$ & 22&23&24&25&26&27&28&29&30&31&32&34&38\\
	\midrule
	prime orbits & 1&1&8&3&7&7&7&3&12&2&3&5&1\\
	\bottomrule
\end{tabular}
\captionof{table}{The $60$ toric degenerations of $E_6/P_6$ by the number of
facets of their polytopes.}
\label{tab:facets}
\end{center}

Three of these degenerations have a known representation-theoretic origin.
Ordering the $27$ weights by $\lambda\le\mu$ when $\mu-\lambda$ is a nonnegative combination of simple roots makes $\pL$ a distributive lattice, whose poset $P_{E_6}$ of join-irreducibles has $16$ elements and $20$ cover relations. 
Standard monomial theory degenerates $E_6/P_6$ flatly to the \emph{Hibi toric variety} of this lattice~\cite{GonciuleaLakshmibai}, whose ideal is generated by the binomials $p_\lambda p_\mu-p_{\lambda\wedge\mu}p_{\lambda\vee\mu}$ for incomparable $\lambda,\mu$, and whose polytope is the order polytope of $P_{E_6}$~\cite{Hibi}. 
Exactly one of our $60$ orbits is this degeneration: up to $W(E_6)$ and signs, its initial ideal is precisely the ideal of the $27$ Hibi relations, and it is the unique orbit whose polytope has the minimal number $22=20+2$ of facets. 
Exactly two further orbits give polytopes unimodularly equivalent to Stanley's chain polytope of $P_{E_6}$~\cite{Stanley}, with $28$ facets. 
We do not know a representation-theoretic description of the remaining $57$ orbits; it would be interesting to realize them as marked poset polytopes or as Newton--Okounkov bodies of explicit valuations.

\section*{Code}
The code and data for the computations in this abstract are available at:
\begin{center}
	\url{https://dcorey2814.github.io/E6P6CoxeterMatroids/}
\end{center}

\section*{AI disclosure}

The author used Gemini (Gemini~3.1), ChatGPT (Sol) and Claude (Opus~5.5) for brainstorming, for finding new proofs of results that the author had obtained independently, for drafting this paper, and for adapting the author's code into the companion GitHub repository and website. Notably, Gemini discovered the tritangent exchange property, Claude formulated the leaf order proof of Theorem \ref{thm:main} and found the representation-theoretic origin of three of the toric degenerations. The author checked all claims made by these tools.

\bibliography{corey_e6p6}

\end{document}